\documentstyle[11pt]{article}
\def\inO #1{\mbox{\bf O\hskip -7pt\raise 1.5pt\hbox{\scriptsize #1}}}
\def\inSO #1{\mbox{\bf SO\hskip -7pt\raise 1.5pt\hbox{\scriptsize #1}}}

\begin{document}
\bibliographystyle{unsrt}

\vbox{\vspace{38mm}}

\begin{center}
{\Large \bf HIGHER ORDER RECURRENCES FOR ANALYTICAL FUNCTIONS OF
TCHEBYSHEFF TYPE
 \\[5mm]}

 A. K. Kwa\'sniewski \\[3mm]
{\it Higher School of Mathematics and Applied Informatics\\
PL - 15-021 Bialystok , ul.Kamienna 17,  Poland\\

e-mail: kwandr @ uwb.edu.pl\\
}

\end{center}

\begin{abstract}
Relation of hyperbolons of volume one to generalized Clifford
algebras is described in [1b] and there some applications are
listed.

\noindent In this note which is an extension of [8] we use the one
parameter subgroups of the group of hyperbolons of volume one in
order to define and investigate generalization of Tchebysheff
polynomial system.

\noindent Parallely functions of roots of polynomials of any
degree are studied as possible generalization of symmetric
functions considered by Eduard Lucas.

\noindent It is found how functions of roots of polynomial of any
degree are related to this generalization of Tchebysheff
polynomials. The relation is explicit.

\noindent In a primary sense the considered generalization is in
passing from $Z_2$ to $Z_n$ group decomposition of the
exponential.

\noindent We end up with an application of the discovered
generalization to quite large class of dynamical systems with
iteration.

\end{abstract}

\section{Introduction}

\subsection*{Part One: ``Tchebycheff {\it m}-polynomials'' recurrence
equation.}

In order to establish a notation we quote the well known definitions:

$$\label{(1.1)}
\begin{array}{l}
T_n\left( x\right) =\cosh \left( n\alpha \right) ;\;\;\cosh \alpha =x; \\
U_n\left( x\right) =\displaystyle{\frac{\sinh \left( n\alpha
\right) }{\sinh \alpha }} ;\;\cosh \alpha =x.
\end{array}
\eqno{(1.1)}%
$$
Of course
$$
\label{(1.2)} \left[ T_n\left( x\right) +\sqrt{x^2-1}U_n\left( x\right)
\right] ^r=T_{nr}\left( x\right) +\sqrt{x^2-1}U_{nr}\left( x\right)
\eqno{(1.2)}%
$$
It is however more useful to use ``a\&b'' notation :
$$
\label{(1.3)} a_x\left( n\right) =T_n\left( x\right) {\rm \ \ and\ \ }%
b_x\left( n\right) =\sqrt{x^2-1}U_n\left( x\right) \eqno{(1.3)}%
$$
convenient also for the introduction of de Moivre matrix group [1b]:
$$
\label{(1.4)} {\bf M}_x=\left\{ \left(
\begin{array}{cc}
a_x\left( n\right) & b_x\left( n\right) \\
b_x\left( n\right) & a_x\left( n\right)
\end{array}
\right) \equiv M_x\left( n\right) \right\} _{n\in Z};\ \det M_x\left(
n\right) =1 \eqno{(1.4)}%
$$
where (see [3]; formula (50))
$$\label{(1.5)}
\begin{array}{l}
a_x\left( -n\right) =a_x\left( n\right) \\
b_x\left( -n\right) =-b_x\left( n\right)
\end{array}
\eqno{(1.5)}%
$$
Naturally $M_x\left( n\right) M_x\left( m\right) =M_x\left( n+m\right) $
hence
$$
\label{(1.6)}
\begin{array}{l}
a_x\left( n+m\right) =a_x\left( n\right) a_x\left( m\right) +b_x\left(
n\right) b_x\left( m\right) \\
a_x\left( n-m\right) =a_x\left( n\right) a_x\left( m\right) -b_x\left(
n\right) b_x\left( m\right)
\end{array}
\eqno{(1.6)}%
$$
and therefore:
$$
\label{(1.7)}
\begin{array}{l}
a_x\left( n+1\right) =2a_x\left( 1\right) a_x\left( n\right) -a_x\left(
n-1\right) \\
a_x\left( 0\right) =1\;\;a_x\left( 1\right) =\cosh \alpha \equiv x
\end{array}
\ . \eqno{(1.7)}%
$$

This is the recurrence relation for $T_n\left( x\right) =a_x\left( n\right) $
Tchebycheff polynomials.

\noindent In this note we shall investigate a following
generalization of the recurrence relation (1.7) for
$T_n\left(x\right) =a_x\left(n\right)$; namely ``Tchebycheff
$m$-polynomials'' (polynomials in {\it m} constrained variables)
are defined by recurrence equation:
$$
\label{(1.8)} xT_n\left( x\right) =\frac 1m\sum_{s\in Z_m}T_{n+\omega
^s}\left( x\right) \ , \eqno{(1.8)}%
$$
where $\omega =\exp \left( \frac{2\pi i}m\right) $, $m\geq 2$, $T_n\left(
x\right) \equiv a_x\left( n\right) $ and $a_x\left( n\right) $ is to be
defined soon with the help of hyperbolic functions of higher order.

In the case of $m=2$, (1.8) coincides with (1.7), if appropriate initial
conditions are added and in order to derive (1.8) one uses de Moivre group
[1b] for the arbitrary $m\geq 2$ in the similar way as above. Hyperbolic $%
\cosh $ and $\sinh $ functions -- and hence Tchebycheff polynomials of the
first and second type -- are known to be strictly related to symmetric
functions of the roots of quadratic equations (see  (5)  in  [3]).  Lucas
has studied in [3] the relationship between these symmetric functions and
the theories of divisibility, continued fractions, combinatorial analysis,
determinants, quadratic diophantine analysis, continued radicals, etc.\\

\begin{itemize}
\item[] His work was also aimed to be {\it ``the starting point
for a more complete study of the properties of the symmetric
functions of the roots of an algebraic equation with rational
coefficients of any degree''} [3].\\
\end{itemize}

\noindent
To do that one uses in this note the properties of hyperbolic functions
of higher order [1a,b],[2].

\subsection*{Part Two: hyperbolic functions of higher order.}

In the following we shall introduce the hyperbolic functions of arbitrary $%
m\geq 2$ order using projection operators [2] which are important objects on
their own.

Let us define this family of projection operators $\left\{ \Delta _k\right\}
_{k\in Z_m}$ acting on the linear space of functions of complex variable
accordingly to:
$$
\label{(1.9)} \Delta _k:=\frac 1m\sum_{s\in Z_m}\omega
^{-ks}\Omega ^s\ ; \eqno{(1.9)}$$ where  \qquad \qquad \qquad
\qquad \qquad $\left( \Omega f\right) \left( z\right) :=f\left( \omega z\right) \ $.\\

\noindent We used  $\omega$  to the negative power because of
historical reasons - positive power shall appear more accurate.

\noindent The set $\left\{\Delta_k\right\}_{k\in Z_m}$ constitutes
the family of orthogonal projection operators~[2]:

$$
\label{(1.10)} \Delta_i \Delta_j=\delta_{ij}\Delta_i\ ;\quad k\in
Z_m\ . \eqno{(1.10)}
$$

With the help of these projection operators $\left\{ \Delta
_k\right\} _{k\in Z_m}$ we define eigenfunctions of the $\Omega $
operator acting on various linear spaces of functions. Here there
are two examples:

\noindent
1. $\left\{ h_k\left( z\right) \right\} _{k\in Z_m}$ constitute the set of
eigenfunctions of the $\Omega $ ($h_k\leftrightarrow $ $m$-hyper\-bolic $k$%
-series);
$$
h_k:=\Delta _k\exp \Rightarrow h_k\left( z\right) =\sum_{s\geq 0}\frac{%
z^{ms+k}}{\left( ms+k\right) !}\Rightarrow \Omega h_k=\omega ^kh_{k\ };\ \
k\in Z_m\ .
$$
\noindent
2. $\left\{ g_k\left( z\right) \right\} _{k\in Z_m}$ constitute the set of
eigenfunctions of the $\Omega $ ($g_k\leftrightarrow $ $m$-geometric $k$%
-series);
$$
g_k:=\Delta _k\frac 1{1-id}\left( {\rm where\ }\frac 1{1-id}\left( z\right)
:=\frac 1{1-z}\right) \Rightarrow g_k\left( z\right) =$$
$$=\sum_{s\geq
0}z^{ms+k}\Rightarrow \Omega g_k=\omega ^kg_{k\ };\ \ k\in Z_m\ .
$$

\noindent
The operators with the algebraic properties of $\Omega $ had been already
profitably used in [4] and [5].

\noindent
The eigenfunctions $\left\{ h_s\left( z\right) \right\} _{s\in Z_m}$ of the $%
\Omega $ ($h_s\leftrightarrow $ $m$-hyperbolic $s$-series) are known since
late 40's as hyperbolic functions of $m$-th order (see [1a,b] for the
references).

\noindent
These generalizations of cosh and sinh hyperbolic functions are given
explicitly by:
$$\label{(1.11)} h_i\left( x\right) =\frac 1m\sum_{k\in {\bf Z}_m}\omega
^{-ki}\exp \left\{ \omega ^kx\right\} ;\;\;i\in {\bf Z}_m;\;\;\omega =\exp
\left\{ i\frac{2\pi }m\right\} \eqno{(1.11)}%
$$
We call (1.11) --- Euler formulae for hyperbolic functions of $m$-th order.
In the next section we construct an analogue of Tchebycheff polynomial
system -- with the help of
$$
\label{(1.12)} h_0\left( x\right) =\frac 1m\sum_{k\in {\bf Z}_m}\exp \left\{
\omega ^kx\right\} ,\quad m\geq 2\ . \eqno{(1.12)}%
$$
\section{Tchebycheff Polynomials' Systems for Higher Order Equations}

One may define Tchebycheff -like functions via recurrences in the same
manner as in (1.7) due to the identity
$$
\label{(2.1)} \sum_{k\in Z_m}h_0\left( \alpha +\omega ^k\beta \right) \equiv
mh_0\left( \alpha \right) h_0\left( \beta \right) . \eqno{(2.1)}%
$$
\indent
We define then the Tchebycheff polynomials' systems of {\it m}-th order to
be given by
$$
\label{(2.2)} T_{\vec n}\left( x\right) =h_0\left( \vec n\alpha \right)
,\quad h_0\left( \alpha \right) \equiv x\ ; \eqno{(2.2)}%
$$
where $\ \qquad \qquad \qquad \vec n\in \left\{
n,n+\omega ,...,n+\omega ^{m-1};\ n\in N\right\} $.\\

The recurrence relation (1.8) for the system $\left\{ T_{\vec n}\left(
x\right) \right\} _{\vec n\in \vec N}$ of Tchebycheff-like functions follows
from (2.1) and (2.2):
$$
xT_n\left( x\right) =\frac 1m\sum_{s\in Z_m}T_{n+\omega ^s}\left( x\right) \
.
$$

As a sufficient example we write this recurrence for the case of $m=3$:
$$
\label{(2.3)} \left\{
\begin{array}{l}
T_{n+1}\left( x\right) =3xT_n\left( x\right) -T_{n+\omega }\left( x\right)
-T_{n+\omega ^2}\left( x\right) \\
T_0\left( x\right) =1\quad T_1\left( x\right) =x
\end{array}
\ .\right. \eqno{(2.3)}%
$$
It is easily seen that $e^\alpha, e^{\omega \alpha}, e^{\omega ^2\alpha }$
are solutions of the characteristic equation of the recurrence (1.8)
$$
\label{(2.4)} \lambda =3x-\lambda ^\omega -\lambda ^{\omega ^2}. \eqno{(2.4)}%
$$
under the identification $x=h_0\left( \alpha \right) $; hence the Binet form
of the solution of the recurrence (1.8) for $m=3$ and under the
identification $x=h_0\left( \alpha \right) $ reads as follows
$$
\label{(2.5)} \frac{\left( e^\alpha \right) ^n+\left( e^{\omega \alpha
}\right) ^n+\left( e^{\omega ^2\alpha }\right) ^n}3\equiv h_0\left( n\alpha
\right) \eqno{(2.5)}%
$$
In order to solve the recurrence (1.8) it is necessary to introduce for $m=3$
:
$$
x=h_0\left( \alpha \right) ;\quad x^{*}=h_0\left( -\alpha \right) ;\quad
x^{**}=h_0\left( \left( 2+\omega \right) \alpha \right) .
$$
Then the following identities and identifications for $m=3$ are useful:

\begin{tabular}{ll}
  0. $T_\beta \left( x\right) :=h_0\left( \beta \alpha \right) ;
    \beta \in {\bf C}$; & \\
  1. $x^{*}=x=x^{**}\quad iff\quad m=2$; &
    2. $h_0\left( \omega ^l\alpha \right)=h_0\left( \alpha \right) ;
    \quad l\in {\bf Z}_3$ ; \\
  3. $T_{\omega ^2}\left( x\right) =T_\omega \left( x\right) =T_1\left(
    x\right) =x$; & 4. $T_{2+\omega }\left( x\right) =x^{**}$; \\
  5. $T_{1+\omega ^2}\left( x\right) =T_{1+\omega }\left( x\right) =x^{*}$; &
    6. $T_{2+\omega }\left( x\right) +T_{2+\omega ^2}\left( x\right)
    =3xx^{*}-1$.
\end{tabular}\\

\noindent
In order to generate these formulae use the identity and the convolution
formula [1b]:
$$
\label{(2.6)} \sum_{k\in Z_3}h_k\left( \alpha \right) h_{i-k}\left( \beta
\right) \equiv h_i\left( \alpha +\beta \right) . \eqno{(2.6)}%
$$

Let us now define ($m=3$) ordinary generating functions for the sequence $%
\left\{ T_{\vec n}\left( x\right) \right\} _{\vec n\in \vec N}$ of
Tchebycheff ``3-polynomials'' of the 0-th kind $(Z_m$ -- group enumeration
terminology) as follows:
$$
\label{(2.7)}
\begin{array}{c}
T_{\left( 0\right) }\left( x,x^{*};z\right) =\displaystyle{\sum_{n\geq
0}}T_n\left(
x,x^{*}\right) z^n\ , \\
T_{\left( s\right) }\left( x,x^{*},x^{**};z\right) =\displaystyle{\sum_{n\geq
0}}T_{n+\omega ^s}\left( x,x^{*},x^{**}\right) z^n\ ,\quad s=1,2\ .
\end{array}
\eqno{(2.7)}%
$$

Using the above (2,..,6) identities and identifications one proves in a
standard way that for the sequence of Tchebycheff ``3-polynomials'' from the
``main stream'' $\left\{ T_n\left( x,x^{*}\right) \right\}
_{n=0}^\infty $ we have
$$
\label{(2.8)} T_{\left( 0\right) }\left( x,x^{*};z\right) =\frac{%
1-2xz+x^{*}z^2}{1-3xz+3x^{*}z^2-z^3}\ ; \eqno{(2.8)}%
$$
while ordinary generating functions for ``aside streams'' are given by
$$
\label{(2.9)} T_{\left( 1\right) }\left( x,x^{*},x^{**};z\right) =\frac{%
x-\left( 3x^2-x^{*}\right) z+x^{**}z^2}{1-3xz+3x^{*}z^2-z^3}\ ; \eqno{(2.9)}%
$$
$$
\label{(2.10)} T_{\left( 2\right) }\left( x,x^{*},x^{**};z\right) =\frac{%
x-\left( 3x^2-x^{*}\right) z+\left( 3xx^{*}-x^{**}-1\right) z^2}{%
1-3xz+3x^{*}z^2-z^3}\ ; \eqno{(2.10)}%
$$
We conclude that for $m>2$ one gets the system of polynomials in {\it m}
interdependent variables.

\section{Cyclic-Symmetric Functions of Roots of Polynomials and
Hyperbolic Functions}

One may show (see formulae (3.6)) that the appropriate analogues of (1.6)
formulae exist (of course $m$ of them instead of two); also the relation of
Tchebysheff $m$-polynomials in $m$ dependent variables to functions of roots
of polynomials generalizes to the case of $m>2$.

\noindent
We are now going to elaborate more on that hyperbolic-trigonometric
character of the introduced Tchebysheff-like systems. To do that, let
us recall de Moivre formulae in their matrix form [1]:
$$
\label{(3.1)}
\begin{array}{l}
H\left( \alpha \right) H\left( \beta \right) =H\left( \alpha +\beta \right)
, \\
H\left( \alpha \right) =\exp \left\{ \gamma \alpha \right\} , \\
\quad {\rm where\ }\gamma =\left( \delta _{i,k\dot -i}\right) ;\ k,i\in Z_m
,
\end{array}
\eqno{(3.1)}%
$$
equivalent to their convolution form [1b]
$$
\label{(3.2)} h_k\left( \alpha +\beta \right) =\sum_{i\in Z_m}h_i\left(
\alpha \right) h_{k\dot -i}\left( \beta \right) ;\quad k\in Z_m\ .
\eqno{(3.2)}%
$$
Then for $m=3$ we have
$$
\label{(3.3)} H\left( n\alpha \right) =\left(
\begin{array}{ccc}
h_0\left( n\alpha \right) & h_1\left( n\alpha \right) & h_2\left( n\alpha
\right) \\
h_2\left( n\alpha \right) & h_0\left( n\alpha \right) & h_1\left( n\alpha
\right) \\
h_1\left( n\alpha \right) & h_2\left( n\alpha \right) & h_0\left( n\alpha
\right)
\end{array}
\right) \equiv H^n\left( \alpha \right) \ . \eqno{(3.3)}%
$$
In order to present the main idea we restrict our attention from
now on to the
case of $m=3$ (the generalization to the arbitrary $m$ is straightforward)
and we introduce the ``$a,b,c$'' notation:
$$
\label{(3.4)} H\left( n\alpha \right) \equiv M_\alpha \left( n\right)
=\left(
\begin{array}{ccc}
a_\alpha \left( n\right) & b_\alpha \left( n\right) & c_\alpha \left(
n\right) \\
c_\alpha \left( n\right) & a_\alpha \left( n\right) & b_\alpha \left(
n\right) \\
b_\alpha \left( n\right) & c_\alpha \left( n\right) & a_\alpha \left(
n\right)
\end{array}
\right)  \eqno{(3.4)}%
$$
thus obtaining de Moivre group [1b]:
$$
\label{(3.5)} M_\alpha =\left\{ M_\alpha \left( n\right) \right\} _{n\in Z}\
. \eqno{(3.5)}%
$$
Due to (3.2) or equivalently due to the group property of (3.5) one obtains (%
$k,n\in Z$):
$$
\label{(3.6a)} a_\alpha \left( n+k\right) =3a_\alpha \left( k\right)
a_\alpha \left( n\right) -a_\alpha \left( n+k\omega \right) -a_\alpha \left(
n+k\omega ^2\right) \eqno{(3.6a)}%
$$
$$\label{(3.6b)} b_\alpha \left( n+k\right) =3a_\alpha \left( k\right)
b_\alpha \left( n\right) -b_\alpha \left( n+k\omega \right) -b_\alpha \left(
n+k\omega ^2\right) \eqno{(3.6b)}%
$$
$$\label{(3.6c)} c_\alpha \left( n+k\right) =3a_\alpha \left( k\right)
c_\alpha \left( n\right) -c_\alpha \left( n+k\omega \right) -c_\alpha \left(
n+k\omega ^2\right) \eqno{(3.6c)}%
$$

i.e. the recurrent sequences $\left\{a_\alpha \left( n\right)
\right\}_{n\in Z}$, $\left\{ b_\alpha \left( n\right) \right\}
_{n\in Z}$, $\left\{ c_\alpha \left( n\right) \right\}_{n\in Z}$
obey the same rule of formation (put $k=1$) but differ in initial
conditions. Hence we have Tchebycheff 3-polynomials of the
zero-th, the first and the second kind $(Z_3$- group enumeration
terminology) where, of course, the recurrence (3.6a) coincides
with (2.3) for $k=1$.

\noindent
Parallely we shall propose a generalization of Lucas symmetric functions [3].

For the convenience of presentation put $m=3$ and consider $V$, $U$, $W$
functions of roots of the equation: (the case $a=b=c$ is excluded)
$$
\label{(3.7)} \left( x-a\right) \left( x-b\right) \left( x-c\right)
=0\Leftrightarrow x^3=Px^2+Qx+R \eqno{(3.7)}%
$$
as follows
$$
\label{(3.8)} V_n\left( a,b,c\right) =a^n+b^n+c^n \eqno{(3.8)}%
$$
where $\qquad \qquad \qquad \qquad \qquad \left\{
\begin{array}{l}
V_{n+3}=PV_{n+2}+QV_{n+1}+RV_n \\
V_0=3,\quad V_1=-P,\quad V_2=a^2+b^2+c^2
\end{array}
\right. ;$
$$
\label{(3.9)} U_n\left( a,b,c\right) =\frac{a^n+\omega b^n+\omega ^2c^n}{%
a+\omega b+\omega ^2c} \eqno{(3.9)}%
$$
where $\qquad \qquad \qquad \left\{
\begin{array}{l}
U_{n+3}=PU_{n+2}+QU_{n+1}+RU_n \\
U_0=0,\quad U_1=1,\quad U_2=\displaystyle{\frac{a^2+\omega b^2+\omega
^2c^2}{a+\omega
b+\omega ^2c}}
\end{array}
\right. ;$
$$
\label{(3.10)} W_n\left( a,b,c\right) =\frac{a^n+\omega ^2b^n+\omega c^n}{%
a+\omega ^2b+\omega c} \eqno{(3.10)}%
$$
where $\qquad  \qquad \qquad \left\{
\begin{array}{l}
W_{n+3}=PW_{n+2}+QW_{n+1}+RW_n \\
W_0=0,\quad W_1=1,\quad W_2=\displaystyle{\frac{a^2+\omega ^2b^2+\omega
c^2}{a+\omega
^2b+\omega c}}
\end{array}
\right. ;$

\medskip
One identifies the values of $P,Q,R$ from (3.7) to be :
$$
\label{(3.11)} -P=a+b+c;\quad Q=ab+ac+bc; \quad -R=abc. \eqno{(3.11)}%
$$

Note that $V_n\left( a,b,c\right) $ -- functions are symmetric, while $%
U_n\left( a,b,c\right) =$\break
$U_n\left( b,c,a\right) =U_n\left( c,a,b\right), $
i.e. $U_n\left( a,b,c\right) $ -- functions are cyclic-symmetric. Similarly $%
W_n\left( a,b,c\right) =W_n\left( b,c,a\right) =W_n\left( c,a,b\right) $.

Recall now that Lucas considered the following symmetric functions of
the roots $a$ and $b$ of the quadratic equation $z^2=Pz-Q$ with $P,Q\in Z$
and $P,Q$ relatively prime.
$$
\begin{array}{l}
V_n\left( a,b\right) =V_n\left( b,a\right) ; \\
U_n\left( a,b\right) =U_n\left( b,a\right) ;
\end{array}
$$
$$
\label{(3.12)}
\begin{array}{l}
V_n\left( a,b\right) =a^n+b^n \\
U_n\left( a,b\right) =\displaystyle{\frac{a^n-b^n}{a-b} }
\end{array}
\eqno{(3.12)}%
$$
Lucas proved [3] the complete analogy of the $V_n$ and $U_n$
symmetric functions of roots with the circular and hyperbolic functions of $%
m=2$ order due to the following ``{\it Lucas formulae}'':
$$
\label{(3.13)} V_n\left( a,b\right) =2Q^{\frac n2}\cosh \left[ \frac n2\ln
\frac ab\right] ; \eqno{(3.13)}%
$$
$$
\label{(3.14)} U_n\left( a,b\right) =\frac{2Q^{\frac n2}}{\sqrt{\Delta }}%
\sinh \left[ \frac n2\ln \frac ab\right] ; \eqno{(3.14)}%
$$
where $\Delta =P^2-4Q$. The equations (3.11) and (3.12) define a one to one
correspondence between formulae in plane (hyperbolic) trigonometry and
analogous formulae for symmetric functions $V_n\left( a,b\right) $ and $%
U_n\left( a,b\right) $. (Below we shall propose an available generalization
of these formulae to the case of arbitrary $m$).

Due to (3.13) and (3.14) we get the following identification (via rescaling
of $V$'s \& $U$'s):
$$
\label{(3.15)} a_\alpha \left( n\right) :=Q^{-\frac n2}\frac
12V_n;\;\;b_\alpha \left( n\right) :=Q^{-\frac n2}\frac{\sqrt{\Delta }}2U_n
\eqno{(3.15)}%
$$
According to this (see also: formulae (49), (50), (51) in [3]) one has de
Moivre group for $m=2$
$$
\label{(3.16)} \left\{ \left(
\begin{array}{cc}
a_\alpha \left( n\right) & b_\alpha \left( n\right) \\
b_\alpha \left( n\right) & a_\alpha \left( n\right)
\end{array}
\right) ;n\in Z\right\} =\left\{ \left(
\begin{array}{cc}
a_\alpha \left( 1\right) & b_\alpha \left( 1\right) \\
b_\alpha \left( 1\right) & a_\alpha \left( 1\right)
\end{array}
\right) ^n;n\in Z\right\} \ ; \eqno{(3.16)}%
$$
where
$$
\label{(3.17)} \alpha =\frac 12\left[ \ln a-\ln b\right] \ . \eqno{(3.17)}%
$$
One may show that it is possible to extend Lucas formulae to the case of the
cyclic-symmetric functions of the roots of the polynomial with rational,
real or complex coefficients of any degree.

For that purpose let us consider cyclic-symmetric functions defined by
(3.8), (3.9), (3.10).

We state that the following identifications hold (generalization to
arbitrary $m$ is trivial):
$$
\label{(3.18a)} a_\alpha \left(n\right)=R^{-\frac n3}
\frac{V_n \left(A,A^\omega,A^{\omega ^2}\right)}3 ; \eqno{(3.18a)}%
$$
$$
\label{(3.18b)} c_\alpha \left( n\right) =R^{-\frac n3}
\frac{U_n\left(
A,A^\omega, A^{\omega ^2}\right) }3h_1\left( \ln A\right) ; \eqno{(3.18b)}%
$$
$$
\label{(3.18c)} b_\alpha \left( n\right) =R^{-\frac n3}\frac{W_n\left(
A,A^\omega, A^{\omega ^2}\right) }3h_2\left( \ln A\right) ; \eqno{(3.18c)}%
$$
Here
$$
\label{(3.19)} \alpha =\frac 13\left[ \ln a+\omega \ln b+\omega ^2\ln
c\right] \ ; \eqno{(3.19)}%
$$
and
$$
\label{(3.20)} A^3=ab^\omega c^{\omega ^2}. \eqno{(3.20)}%
$$
Accordingly to this one has de Moivre group for $m=3$:
$$
\label{(3.21)} H\left( n\alpha \right) \equiv \left(
\begin{array}{ccc}
a_\alpha \left( n\right) & b_\alpha \left( n\right) & c_\alpha \left(
n\right) \\
c_\alpha \left( n\right) & a_\alpha \left( n\right) & b_\alpha \left(
n\right) \\
b_\alpha \left( n\right) & c_\alpha \left( n\right) & a_\alpha \left(
n\right)
\end{array}
\right) \ ;\qquad n\in Z \eqno{(3.21)}%
$$
The corresponding formulae for the case of $m=2 (\omega =-1,
h_1\equiv \sinh  )$ are:
$$
\label{(3.18*a)} a_\alpha \left( n\right) =Q^{-\frac n2}\frac{V_n\left(
A,A^\omega \right) }2\ ; \eqno{(3.18*a)}%
$$
$$
\label{(3.18*b)} b_\alpha \left( n\right) =Q^{-\frac n2}\frac{U_n\left(
A,A^\omega \right) }2h_1\left( \ln A\right) \ ; \eqno{(3.18*b)}%
$$
with
$$
\label{(3.19*)} \alpha =\frac 12\left[ \ln a+\omega \ln b\right] \ ;
\eqno{(3.19*)}%
$$
and
$$
\label{(3.20*)} A^2=ab^\omega \ . \eqno{(3.20*)}%
$$
Of course for the $Q=1$ case (via obvious rescaling of $V$'s \& $U$'s one
may always get to this case):
$$
\label{(3.22)} V_n\left( a,b\right) =V_n\left( A,A^\omega \right) \ ;
\eqno{(3.22)}%
$$
$$
\label{(3.23)} U_n\left( a,b\right) =U_n\left( A,A^\omega \right) \ .
\eqno{(3.23)}%
$$
and both $\left\{ a_\alpha \left( n\right) \right\} _{n\in Z}${\bf,}
$\left\{ b_\alpha \left( n\right) \right\} _{n\in Z}$ sequences and $\left\{
V_n\left( a,b\right) \right\} _{n\in Z}$,\break $\left\{ U_n\left( a,b\right)
\right\} _{n\in Z}$ recurrent sequences obey the same rule of formation but
differ in initial conditions.  For $m=3$ we have
$$
\label{(3.24a)} V_n\left( A,A^\omega, A^{\omega ^2}\right) \neq V_n\left(
a,b,c\right) \ , \eqno{(3.24a)}%
$$
$$
\label{(3.24b)} U_n\left( A,A^\omega, A^{\omega ^2}\right) \neq U_n\left(
a,b,c\right) \ , \eqno{(3.24b)}%
$$
$$
\label{(3.24c)} W_n\left( A,A^\omega, A^{\omega ^2}\right) \neq W_n\left(
a,b,c\right) \ . \eqno{(3.24c)}%
$$

Hence for $m=3$ the ``proper'' generalization of Lucas formulae (3.13) and
(3.14) or equivalently (3.15) is given by formulae (3.18).

\noindent
Functions $V_n\left( A,A^\omega, A^{\omega ^2}\right) $ are of course
symmetric in $A, A^\omega, A^{\omega ^2}$ arguments but are no more symmetric
functions of roots $a,b,c$ ; they are however cyclic symmetric functions of
roots.

\noindent
Functions $U_n\left( A,A^\omega, A^{\omega ^2}\right) $, $W_n\left(
A,A^\omega, A^{\omega ^2}\right) $ are of course cyclic-symmetric in $%
A,A^\omega, A^{\omega ^2}$ arguments, but are no more cyclic- symmetric
functions of roots $a,b,c$.

\section{Geometrical Representation of Tchebycheff {\it m}-Polynomials
of {\it k}-th Kind.}

Consider the case $m=2$\ \ i.e. the ordinary Tchebycheff polynomials in one
variable $x=\cosh \alpha $ or in two dependent variables: $x=\cosh \alpha $
and $y=\sinh \alpha $ then one has
$$
\label{(4.1)} T_n\left( x,y\right) =\sum_{k=0}^{\left[ \frac n2\right]
}\left(
\begin{array}{c}
n \\
2k
\end{array}
\right) x^{n-2k}y^{2k}, \eqno{(4.1)}%
$$
where $x$, $y$ are coordinates of a point from a hyperbola given by the
group of hyperbolons of volume one [1b]
$$
\label{(4.2)} \det \left(
\begin{array}{cc}
x & y \\
y & x
\end{array}
\right) =1\ . \eqno{(4.2)}%
$$

Consider now the case $m=3$. Then one may show that
$$
\label{(4.3)} h_0\left( nx\right) =\sum_{k=0}^n\left(
\begin{array}{c}
n \\
k
\end{array}
\right) \sum_{i=0}^{n-k}\left(
\begin{array}{c}
n-k \\
i
\end{array}
\right) \delta \left( i\dot -k\right) h_0^{n-k-i}\left( x\right) h_1^i\left(
x\right) h_2^k\left( x\right) , \eqno{(4.3)}%
$$
i.e. for Tchebycheff 3-polynomials of zero-th kind we get
$$
\label{(4.4)} T_n\left( x,y,z\right) \equiv \sum_{k=0}^n\left(
\begin{array}{c}
n \\
k
\end{array}
\right) \sum_{i=0}^{n-k}\left(
\begin{array}{c}
n-k \\
i
\end{array}
\right) \delta \left( i\dot -k\right) x^{n-k-i}y^iz^k, \eqno{(4.4)}%
$$
where $x,y,z$ are coordinates of a point from the surface given by the group
of hyperbolons of volume one [1b]
$$
\label{(4.5)} \det \left(
\begin{array}{ccc}
h_0\left( \alpha \right) & h_1\left( \alpha \right) & h_2\left( \alpha
\right) \\
h_2\left( \alpha \right) & h_0\left( \alpha \right) & h_1\left( \alpha
\right) \\
h_1\left( \alpha \right) & h_2\left( \alpha \right) & h_0\left( \alpha
\right)
\end{array}
\right) \equiv \det \left(
\begin{array}{ccc}
x & y & z \\
z & x & y \\
y & z & x
\end{array}
\right) =1. \eqno{(4.5)}%
$$
This surface defined by the equation $x^3+y^3+z^3-3xyz=1$.

\noindent
Tchebycheff 3-polynomials of the first and second kind are also functions on
the one-parameter subgroup of hyperbolons of volume one [1] as seen from the
formulae
$$\label{(4.6)} U_n\left( x,y,z\right) \equiv \sum_{k=0}^n\left(
\begin{array}{c}
n \\
k
\end{array}
\right) \sum_{i=0}^{n-k}\left(
\begin{array}{c}
n-k \\
i
\end{array}
\right) \delta \left( i+1\dot -k\right) x^{n-k-i}y^iz^k, \eqno{(4.6)}%
$$
$$
\label{(4.7)} W_n\left( x,y,z\right) \equiv \sum_{k=0}^n\left(
\begin{array}{c}
n \\
k
\end{array}
\right) \sum_{i=0}^{n-k}\left(
\begin{array}{c}
n-k \\
i
\end{array}
\right) \delta \left( i\dot -\left( k+1\right) \right) x^{n-k-i}y^iz^k.
\eqno{(4.7)}%
$$

The formulae (4.4), (4.6) and (4.7) may be obtained with the help of the
following identities:
$$
\label{(4.8)} \left( a+b+c\right) ^n\equiv \sum_{k=0}^n\left(
\begin{array}{c}
n \\
k
\end{array}
\right) \left( a+b\right) ^{n-k}c^k\equiv \sum_{k=0}^n\left(
\begin{array}{c}
n \\
k
\end{array}
\right) \sum_{i=0}^{n-k}\left(
\begin{array}{c}
n-k \\
i
\end{array}
\right) a^{n-k-i}b^ic^k \eqno{(4.8)}%
$$
$$
\label{(4.9)} \frac 13\left\{ \omega ^0+\omega ^{k+2i}+\omega
^{i+2k}\right\} =\delta \left( i\stackrel{\cdot }{-}k\right) \quad \omega
=\exp \left\{ \frac{2\pi i}3\right\} \eqno{(4.9)}%
$$
where $\qquad \qquad \stackrel{\cdot}{-}\equiv
substracting\  {\rm mod }\ 3$.\\

Tchebycheff $m$-polynomials of the $k$-th kind are of course also functions
on the one-parameter subgroup of hyperbolons of volume one [1b]
in the\break  $m$-th dimensional space -- the whole group being
represented by points of the corresponding surface. For example
in $m=4$ case the group of hyperbolons of volume one is
represented by points of the surface defined by equation [6]
$$
\label{(4.10)}
-x^4+y^4-z^4+t^4+4x^2yt-4xy^2z+4z^2yt-4t^2xz+2x^2z^2-2y^2t^2=1. \eqno{(4.10)}%
$$
\section{Cayley-Hamilton Theorem and Matrix Recurrences for Tchebycheff
{\it m}-Polynomials}

Let us consider the matrix
$$
\label{(5.1)} A=\left(
\begin{array}{cc}
P & -Q \\
1 & 0
\end{array}
\right) . \eqno{(5.1)}%
$$
This matrix $A$ is the representation of equivalence class of all matrices ($%
2\times 2$) with trace $P$ and determinant $Q$. According to Cayley-Hamilton
Theorem it satisfies the equation
$$
\label{(5.2)} A^2=PA-QI\ , \eqno{(5.2)}%
$$
which is of the form $z^2=Pz-Q$\ \ i.e. of the form of characteristic equation
for the corresponding recurrences (compare with (3.12) and with (10) in [3]):
$$
\label{(5.3)} V_{n+2}=PV_{n+1}-QV_n\ ,\quad U_{n+2}=PU_{n+1}-QU_n\ .
\eqno{(5.3)}%
$$
Matrix $A$ generates solutions of these recurrences due to
$$
\label{(5.4)} \left(
\begin{array}{c}
F_{n+2} \\
F_{n+1}
\end{array}
\right) =\left(
\begin{array}{cc}
P & -Q \\
1 & 0
\end{array}
\right) \left(
\begin{array}{c}
F_{n+1} \\
F_n
\end{array}
\right) . \eqno{(5.4)}%
$$
For example if we take $F_0=0,\ F_1=1$ as initial values then [7]
$$
\label{(5.5)} A^n=\left(
\begin{array}{cc}
F_{n+2} & -QF_{n+1} \\ \\
F_{n+1} & -QF_n
\end{array}
\right) ;\qquad n\in N\ . \eqno{(5.5)}%
$$

Let us consider now the matrix
$$
\label{(5.6)} A=\left(
\begin{array}{ccc}
P & Q & R \\
1 & 0 & 0 \\
0 & 1 & 0
\end{array}
\right) . \eqno{(5.6)}%
$$
This matrix $A$ is the representation of equivalence class of all matrices ($%
3\times 3$) with trace $P$, determinant $R$ and sum of corresponding minors $%
=-Q$. According to Cayley-Hamilton Theorem $A$ satisfies
$$
\label{(5.7)} A^3=PA^2+QA+RI\ , \eqno{(5.7)}%
$$
which is of the form $\left( x-a\right) \left( x-b\right) \left( x-c\right)
=0\quad \Leftrightarrow \quad x^3=Px^2+Qx+R.$ i.e. of the form of
characteristic equation for the corresponding recurrences (compare with
(3.8), (3.9) and (3.10)):
$$
\label{(5.8)} \left\{
\begin{array}{l}
F_{n+3}=PF_{n+2}+QF_{n+1}+RF_n \\
initial...values
\end{array}
\right. . \eqno{(5.8)}%
$$
Matrix $A$ generates solutions of these recurrences due to
$$
\label{(5.9)} \left(
\begin{array}{c}
F_{n+3} \\
F_{n+2} \\
F_{n+1}
\end{array}
\right) =\left(
\begin{array}{ccc}
P & Q & R \\
1 & 0 & 0 \\
0 & 1 & 0
\end{array}
\right) \left(
\begin{array}{c}
F_{n+2} \\
F_{n+1} \\
F_n
\end{array}
\right) ;\quad n=0,1,2,...\ . \eqno{(5.9)}%
$$
For example if we take $F_0=0,\ F_1=1,F_2=1$ as initial values then [7]
$$
\label{(5.10)} A^n=\left(
\begin{array}{ccc}
F_{n+2} & QF_{n+1}+RF_n & RF_{n+1} \\
F_{n+1} & QF_n+RF_{n-1} & RF_n \\
F_n & QF_{n-1}+RF_{n-2} & RF_{n-1}
\end{array}
\right) \quad n=2,3,...\ . \eqno{(5.10)}%
$$

In general case of arbitrary $m>1$ one obtains the matrix $A$ adjoint
representation of an endomorphism $A$ -- given by its invariants from
Cayley-Hamilton theorem -- after choosing the basis
$$
\label{(5.11)} \left\{ e_k=A^{m-1-k};\ k\in Z_m\right\} \ . \eqno{(5.11)}%
$$
Of course
$$
\label{(5.12)} A^m=\sum_{k\in Z_m}\alpha _kA^k\ ; \eqno{(5.12)}%
$$
therefore in this general case

$$
\label{(5.13)} A=\left(
\begin{array}{cccccc}
\alpha _{n-1} & \alpha _{n-2} & \alpha _{n-3} & \cdots & \alpha _1 & \alpha
_0 \\
1 & 0 & 0 & \cdots & 0 & 0 \\
0 & 1 & 0 & \cdots & 0 & 0 \\
\vdots & \vdots & \vdots & \ddots & \vdots & \vdots \\
0 & 0 & 0 & \cdots & 1 & 0
\end{array}
\right) . \eqno{(5.13)}%
$$

Matrix $A$ may be considered as a generator of a dynamical system\break
$\left\{
A^n;\ n\in Z\right\} $ in which dynamics is introduced via iteration [7] and
thus is governed by recurrent sequences.

Naturally the appropriate choice of invariants $\left\{ \alpha _k;\ k\in
Z_m\right\} $ leads to Tche\-bycheff $m$-polynomials as exemplified by $m=2$
and $m=3$ cases.

\noindent
This paper extends primary results from [8] -- however it is self-contained.
Further development of the investigation presented above is to be found in
[9] soon.


\begin{thebibliography}{1.b}
\bibitem[1a]{1.a}  Fleury N., M. Rauch de Trautenberg, R. M. Yamaleev,
CRN-PHTH/91-07 (1991) Strasbourg; Universite Louis Pasteur -- preprint.

\bibitem[1.b]{1.b}  Kwa{\'s}niewski A. K., On Hyperbolic and Elliptic
Mappings and Quasi-numbers Algebras, {\it Advances in Applied Clifford
Algebras}, {\bf 2} 107-144, (1992).

\bibitem[2]{2} Kwa{\'s}niewski A. K., On the Onsager Problem for Potts
Models, {\it J. Phys. A: Math. Gen.}, {\bf 19} 1469-1476, (1986).

\bibitem[3]{3}  Lucas Eduard, Th\'eorie des Fonctions Num\'eriques
Simplement P\'eriodiques, {\it American Journal of Mathematics}, {\bf 1}
184-240, (1878); (Translated from the French version by Sidney Kravitz, Edited
by Douglas Lind Fibonacci Association 1969.

\bibitem[4]{4}  Weyl H., ``Theory of groups and Quantum Mechanics'' New
York: E. P. Dutton Co., 1932.

\bibitem[5]{5} Schwinger J., Unitary Operator Basis,
{\it Proc. Nat. Acad. Sci.}, {\bf 46}, 570 (1960).

\bibitem[6]{6}  Yamaleev  R. M., ``Introduction into Theory of N-Unitary
Group'', Communication of the Joint Institute for Nuclear Research,
P2-90-129, Dubna 1990.

\bibitem[7]{7}  Bajguz W. and A. K. Kwa{\'s}niewski, Clifford Algebras,
Quantum Mechanics and Fibonacci-likeSequences, {\it Advances in
Applied Clifford Algebras}, {\bf 4} (1), 73-88, (1994).

\bibitem[8]{8}  Bajguz W. and A. K. Kwa{\'s}niewski,  ``On Hyperbolic
Quasinumbers and Chebyshev-like Functions'' Reports on
Mathematical Physic s{\it 43}(1999): 367-376.

\bibitem[9]{9}   Bajguz W.,  ``On generalization of Tchebycheff
Polynomials'',, Integral Transforms and Special Functions, {\bf 9}
(2), 73-88, (2000):91-98


\end{thebibliography}
\end{document}